\documentclass[12pt,twoside]{article}
\title{On $S$-coherence}
\date{}
\usepackage{amsfonts}
\usepackage{amsmath}
\usepackage{amssymb}
\usepackage{latexsym}
\usepackage{graphicx}
\usepackage[all]{xy}

\newtheorem{thm}{\bf Theorem}[section]
\newtheorem{cor}[thm]{\bf Corollary}
\newtheorem{lem}[thm]{\bf Lemma}
\newtheorem{prop}[thm]{\bf Proposition}
\newtheorem{defn}[thm]{\bf Definition}

\newtheorem{rem}[thm]{\bf Remark}

\newtheorem{exmp}[thm]{\bf Example}

\catcode`\ç=13
\defç{\c{c}}
\catcode`\é=13
\defé{\'e}
\catcode`\à=13
\defà{\`a}
\catcode`\è=13
\defè{\`e}
\catcode`\â=13
\defâ{\^a}
\catcode`\ù=13
\defù{\`u}
\catcode`\ê=13
\defê{\^e}
\catcode`\î=13
\defî{\^\i}
\catcode`\ô=13
\defô{\^o}
\newcommand{\field}[1]{\mathbb{#1}}

\newcommand{\Z }{\field{Z}}
\newcommand{\N }{\field{N}}

\def\proof{{\parindent0pt {\bf Proof.\ }}}
\newcommand{\cqfd}
{\hspace{1cm}
\rule{2mm}{2mm}%
\medbreak%
\par%
}

\def\Im{{\rm Im}}

\author{}
\begin{document}

\thispagestyle{empty}

\maketitle \vspace*{-1.5cm}
\begin{center}{\large\bf Driss   Bennis$^{1,a}$ and Mohammed El
Hajoui$^{1,b}$}

\bigskip

 \small{1. Department of Mathematics, Laboratory of Analysis, Algebra and Decision Support,  Faculty of Sciences,  B.P. 1014,\\
 Mohammed V  University in Rabat,   Rabat, Morocco\\
$\mathbf{a.}$ d.bennis@fsr.ac.ma; driss$\_$bennis@hotmail.com\\
$\mathbf{b.}$  hajoui4@yahoo.fr}
\end{center}

\bigskip

\noindent{\large\bf Abstract.} 
  Recentely,  Anderson and Dumitrescu's $S$-finiteness has attracted the interest of several authors. In this paper, we introduce the notions of $S$-finitely presented modules and then of $S$-coherent rings which are   $S$-versions of  finitely presented modules and  coherent rings, respectively. Among other results, we give an $S$-version of the classical Chase's characterization of coherent rings.  We end the paper with a brief discussion on other $S$-versions of  finitely presented modules and coherent rings. We prove that these last  $S$-versions can be characterized  in terms of localization.\bigskip

\medskip

\bigskip\bigskip


\section{Introduction}   Throughout this paper all rings
are commutative with   identity; in particular, $R$
denotes such a ring, and all modules are unitary.     $S$ will be  a   multiplicative subset of $R$. We use $(I:a)$, for an ideal $I$ and an element $a\in R$, to denote the quotient ideal $\{x\in R;\, xa\in I\}$. \bigskip

According to \cite{AD02}, an $R$ module $M$ is called $S$-finite if there exists a finitely generated submodule $N$ of $M$ such that $sM \subseteq N$  for some $s \in S$. Also, from \cite{AD02}, an $R$-module $M$ is called
$S$-Noetherian if each submodule of $M$ is $S$-finite. In particular,   $R$ is said to be an $S$-Noetherian ring, if it is $S$-Noetherian as an $R$-module; that is, every ideal of $R$ is $S$-finite. It is clear that every Noetherian ring is $S$-Noetherian.\bigskip

The notions of $S$-finite modules and of $S$-Noetherian rings were introduced by Anderson and Dumitrescu motivated by the works done in \cite{HHJ} and  \cite{AKZ95}. They succeeded to generalize several well-known results on Noetherian rings including the classical Cohen's result and Hilbert basis theorem under an additional condition. Since then the $S$-finiteness   has attracted the interest of several authors (see for instance \cite{HH16, HH15, KKL14, LO14, LO15,L07}). Recentely, motivated by the work of Anderson and Dumitrescu,   $S$-versions of some classical notions have been introduced   (see for instance \cite{HH16,KKL14}). In this paper we are inerested in   $S$-versions of finitely presented modules and  coherent rings. Actually, there are two possibilities which could be considered as $S$-versions of finitely presented modules which lead to two $S$-versions of coherent rings.   We prove that the $S$-version of coherent rings defined by one of them  has a characterization similar to the classical one  given by Chase for coherent rings \cite[Theorem  2.2]{C60}. This is why we adopt this notion as the suitable  $S$-version of finitely presented modules. However, it seems not evident to   characterize this notion in terms of localization. We prove that indeed it is the other $S$-version, which is briefly studied at the end of the paper, has a characterization  in terms of localization.\bigskip

The organization of the paper is as follows: In Section 2, we introduce and study an $S$-version of finitely presented modules. We call it an  \textit{$S$-finitely presented module} (see Definition \ref{defn-S-fp}). Then, we study the behavior of $S$-finiteness  in short exact sequences (see Theorem \ref{thm-f-p-sq}). We end Section 2 with some change of rings results (see Proposition \ref{pro-chan-phi} and Corollary \ref{cor-Sfp-Quot}). Section 3 is devoted to the $S$-version of coherent rings which are called \textit{$S$-coherent rings} (see Definition \ref{def-S-C-ring}). Our main result  represents  the $S$-counterpart of Chase's result \cite[Theorem  2.2]{C60} (see Theorem \ref{thm-princ}).   Also an   $S$-version of coherent modules is introduced (see Definition \ref{def-S-C-mod} and Proposition \ref{pro-S-C-sq}). We end the paper with a short section  which  presents the other   $S$-version of $S$-finiteness (see Definitions \ref{def-c-S-fp} and \ref{def-c-S-C-ring}). We prove that these notions can be characterized in terms of localization (see Proposition \ref{prop-c-Sfp-loc} and  Theorem \ref{thm-princ2}). We end the paper with  results which relate   $S$-finiteness    with the notion of $S$-saturation (see Propositions \ref{prop-Sat-S-fp} and \ref{prop-Sat-S-f} and Corollary \ref{cor-princ2}).\medskip


\section{$S$-finitely presented modules}

In this section, we introduce and investigate an $S$-version of the classical finitely presented modules. Other version is discuted in Section 4.

\begin{defn}\label{defn-S-fp} \textnormal{An $R$-module $M$ is called  $S$-finitely presented, if there exists an exact sequence of $R$-modules
 $0\longrightarrow K\longrightarrow F\longrightarrow M\longrightarrow 0$, where $K$ is $S$-finite and $F$ is a finitely generated   free $R$-modules.}
\end{defn}

Clearly, every finitely presented module is $S$-finitely presented. However, the converse does not hold in general. For that,  it suffices to note that when $R$ is a non-Noetherian $S$-Noetherian ring, then  there is an $S$-finite ideal $I$ which is not finitely generated. Then, the $R$-module $R/I$ is  $S$-finitely presented but it is not finitely presented.\\
Also, it is evident that every $S$-finitely presented module is finitely generated. To give an example of a  finitely generated module which is not $S$-finitely presented, it suffices to consider  an ideal $I$ which is not   $S$-finite and then use Proposition   \ref{prop-S-fp-inde} given hereinafter.\medskip

One could remark that in Definition \ref{defn-S-fp} we assume that the free module $F$ is  finitely generated rather than   $S$-finite. In fact, because of the following result both of notions coincide for free modules.

\begin{prop}\label{prop-S-f-free} Every $S$-finite free $R$-module is finitely generated.
\end{prop}
 \proof  Let $M=\bigoplus\limits_{i\in I}Re_i$ be  an  $S$-finite free $R$-module, where $(e_i)_{i\in I}$ is a basis of $M$ and   $I$ is an index set. Then, there exist   a  finitely generated $R$-module $N$ and an $s\in S$ such that $sM\subseteq
N\subseteq M$. Then, $N=Rm_1+\cdots +Rm_n$ for some $m_1,...,m_n\in M$ ($n>0$ is an integer). For  every $k\in \{1,...,n\}$, there exists a  finite subset $J_k$ of $I$ such that $m_k=\sum\limits_{j\in J_k}\lambda_{kj}e_j$.  Let   
$J=\bigcup\limits_{k=1}^n J_k$. Then, the finitely generated $R$-module $M'=\bigoplus\limits_{j\in J}Re_j$ contains  $N$. We show that $M'=M$. Deny. There exists an $i_0\in I \backslash J$ such that
$e_{i_0}\notin M'$. But $se_{i_0}\in N\subseteq M'$ and so  
$se_{i_0}=\sum\limits_{j\in J}\lambda'_{j}e_j$ for some
$\lambda'_j\in R$.    This is impossible since $(e_i)_{i\in I}$ is a basis.\cqfd

\begin{rem}\label{rem-S-fp-proj}  \textnormal{Similarly to the proof of Proposition \ref{prop-S-f-free} above, one can prove that any  $S$-finite torsion-free module cannot be decomposed into an infinite direct sum of non-zero modules. This shows that any $S$-finite projective module is   countably generated by  Kaplansky \cite[Theorem 1]{K58}. Then, naturaly one would ask of the existence of $S$-finite projective module which is not finitely generated. For this, consider the  Boolean ring $R=\prod\limits_{i=1}^{\infty} k_i$, where $k_i$ is the field of two elements for every $i\in \N$. Consider the projective ideal $M=\bigoplus \limits_{i=1}^{\infty} k_i$ and the element $e=(1,0,0,...)$ (see \cite[Example 2.7]{Co94}). Then, $S=\{1,e\}$ is a   multiplicative subset of $R$. Since $eM=k_1$ is a finitely generated $R$-module, $M$ is the desired example of $S$-finite projective module which is not finitely generated.\\
However, determining rings over which every $S$-finite projective module is  finitely generated could be of interest. It is worth noting that rings over which every projective module is a direct sum of finitely generated modules satisfy this condition. These rings were investigated in \cite{MPR7}.}
\end{rem}

Next result shows that, as in the classical case \cite[Lemma 2.1.1]{Glaz}, an $S$-finitely presented module does not depend  on  one specific short exact sequence of the form given  in Definition \ref{defn-S-fp}.

\begin{prop}\label{prop-S-fp-inde} An $R$-module $M$ is $S$-finitely presented if and only of $M$   is  finitely generated and, for every surjective homomorphism of $R$-modules  $F\stackrel{f}{\longrightarrow}M\longrightarrow 0$, where  $F$ is  a finitely generated free $R$-module, $ \ker f $ is $S$-finite.
\end{prop}
 \proof  ($\Leftarrow$) Obvious.\\
($\Rightarrow$)  Since $M$ is  $S$-finitely presented,  
there exists an exact sequence of $R$-modules $0\longrightarrow
K\longrightarrow F'\longrightarrow M\longrightarrow 0$, where $K$ is
$S$-finite and $F'$ is   finitely generated and  free. Then, by Schanuel's lemma,  $K\oplus F\cong\ker f\oplus F'$, then $\ker f$ is $S$-finite.\cqfd

The following result  represnts the behavior of $S$-finiteness in short exact sequences.  It is a generalization  of \cite[Theorem 2.1.2]{Glaz} for modules with $\lambda$-dimension at most $1$. Note that one can give  an $S$-version of the classical $\lambda$-dimension  (see \cite[page]{Glaz}). However, here we prefer to focus on the notion of  $S$-finitely presented modules, and a  discussion on the suitable $S$-version of the $\lambda$-dimension could be the subject of a further work.

\begin{thm}\label{thm-f-p-sq}  Let $0\longrightarrow
M'\stackrel{f}{\longrightarrow} M\stackrel{g}{\longrightarrow}
M''\longrightarrow 0$ be an exact sequence of $R$-modules. The 
following assertions hold:
\begin{enumerate}
\item If $M'$ and $M''$ are
$S$-finite, then $M$ is $S$-finite.\\
 In particular, every finite direct sum of   $S$-finite modules is $S$-finite.
\item If $M'$ and $M''$ are $S$-finitely presented, then $M$ is $S$-finitely presented.\\
In particular, every finite direct sum of $S$-finitely presented modules is $S$-finitely presented.
\item If $M$ is $S$-finite, then $M''$ is $S$-finite.\\
In particular, a direct summand of an $S$-finite module is $S$-finite
\item If $M'$ is $S$-finite and $M$ is $S$-finitely presented, then $M''$ is $S$-finitely presented.
\item If $M''$ is $S$-finitely presented and $M$ is $S$-finite, then $M'$ is $S$-finite.
\end{enumerate}
\end{thm}
\proof 1. Since $M''$ is $S$-finite, there exist  a finitely generated submodule $N''$ of
$M''$ and  an $s\in S$ such that $sM''\subseteq N''$. Let $N''=\sum\limits_{i=1}^nRe_i$  for some $e_i\in M''$ and $n\in \N$. Since $g$ is surjective, there exists an $m_i\in M$ such that $g(m_i)=e_i$ for every $i\in \{1,...,n,\}$. Let $x\in M$, so $sx\in N=g^{-1}(N'')$. Then $g(sx)\in g(N)=N''$, and so $g(sx)=\sum\limits_{i=1}^n\alpha_i e_i=\sum\limits_{i=1}^n\alpha_i g(m_i)=g(\sum\limits_{i=1}^n\alpha_im_i)$. Then, $g(sx-\sum\limits_{i=1}^n\alpha_im_i)=0$. Thus, 
$(sx-\sum\limits_{i=1}^n\alpha_i m_i)\in\ker g=\Im f$ which is $S$-finite. So there exist  a finitely generated submodule  $N'$ of $\Im f$   and an $s'\in S$ such that $s'\Im
f\subseteq N'$. Then, $s'sx\in N'+\sum\limits_{i=1}^nRm_i$ and so $s's M$ is a submodule of $ N'+\sum\limits_{i=1}^nRm_i$ which is a finitely generated submodule of $M$. Therefore, 
$M$ is $S$-finite.\\ 
2. Since $M'$ and $M''$ are $S$-finitely presented, there exist  two shorts
exacts sequences:\\ $0\longrightarrow K'\longrightarrow
F'\longrightarrow M'\longrightarrow0$ and $0\longrightarrow
K''\longrightarrow F''\longrightarrow M''\longrightarrow0$, with
$K'$ and  $K''$ are $S$-finite $R$-modules and $F'$ and $F''$ are  finitely generated  free $R$-modules. Then, by Horseshoe Lemma, we  get the following diagram
$$\xymatrix{&0&0&0&\\0\ar[r]&M'\ar[r]\ar[u]&M\ar[u]\ar[r]&M''\ar[u]\ar[r]
&0\\0\ar[r]&F'\ar[u]\ar@{.>}[r]&F'\oplus
F''\ar@{.>}[u]\ar@{.>}[r]&F''\ar[u]\ar[r]&0\\0\ar[r]&K'\ar[u]\ar@{.>}[r]&K\ar[u]\ar@{.>}[r]&
K''\ar[r]\ar[u]&0\\&0\ar[u]&0\ar[u]&0\ar[u]&}$$ 
By the first assertion,  $K$ is $S$-finite.  Therefore,   $M$ is $S$-finitely presented.\\ 
3. Obvious.\\ 
4. Since $M$ is $S$-finitely presented, there exists a short exact sequence of $R$-modules $0\longrightarrow
K\longrightarrow F\longrightarrow M\longrightarrow0$, where  $K$ is $S$-finite
and $F$ is a finitely generated free $R$-module. Consider the following pullback diagram
$$\xymatrix{&0&0&&\\0\ar[r]&M'\ar[r]\ar[u]&M\ar[u]\ar[r]&M''\ar[r]
&0\\0\ar[r]&D\ar@{.>}[r]\ar@{.>}[u]&F\ar[u]\ar[r]&M''\ar@{=}[u]\ar[r]
&0\\&K\ar[u]\ar@{=}[r]&K\ar[u]&&\\&0\ar[u]&0\ar[u]&&}$$ By (1), $D$ is $S$-finite. Therefore, $M''$ is
$S$-finitely presented.\\ 
5. Since $M''$ is $S$-finitely presented, there exists a short exact sequence $0\longrightarrow K\longrightarrow F\longrightarrow
M''\longrightarrow0$ where $K$ is $S$-finite
and $F$ is a finitely generated free $R$-module.  Consider the following pullback diagram
$$\xymatrix{&&0&0&\\0\ar[r]&M'\ar[r]\ar@{=}[d]&M\ar[u]\ar[r]&M''\ar[u]\ar[r]
&0\\0\ar[r]&M'\ar[r]&D\ar@{.>}[u]\ar@{.>}[r]&F\ar[u]\ar[r]&0\\&&K\ar[u]\ar@{=}[r]&K\ar[u]&&\\&&0\ar[u]&0\ar[u]&}$$
Since $F$ is free,   $D\cong M'\oplus F$, and so $D$ is $S$-finite (since $M'$ and $F$ are $S$-finite). Therefore, $M'$ is $S$-finite.\cqfd

As a simple consequence, we get the following result which extends \cite[Corollary 2.1.3]{Glaz}.

\begin{cor} Let $N_1$ and $N_2$ be two $S$-finitely presented
submodules of an $R$-module. Then,  $N_1+N_2$ is $S$-finitely presented if only if $N_1\cap N_2$ is $S$-finite.
\end{cor}
\proof Use the short exact sequence of $R$-modules $0\longrightarrow N_1\cap N_2\longrightarrow N_1\oplus
N_2\longrightarrow N_1+N_2\longrightarrow 0$.\cqfd

We end this section with the following change of rings results.\medskip

The following result extends \cite[Theorem 2.1.7]{Glaz}.

\begin{prop}\label{pro-chan-phi} Let $A$ and $B$ be
rings, let $\phi:A\longrightarrow B$ be a ring homomorphism making $B$ a  finitely generated $A$-module and let  $V$ be a  multiplicative subset of $A$ such that $0\not\in\phi(V)$. Every $B$-module which is $V$-finitely presented as an $A$-module it is $\phi(V)$-finitely presented as a $B$-module.
\end{prop}
\proof Let $M$ be a $B$-module which is $V$-finitely presented as an $A$-module. Then $M$ is  a finitely generated $A$-module. Then, $M$  is a finitely generated $B$-module. Thus there is an exact sequence of $B$-modules 
$0\longrightarrow K \longrightarrow B^n \longrightarrow  M\longrightarrow0$, where  $n>0$
is an integer.  This sequence   is also  an exact sequence of $A$-modules.  Since $M$ is an   $V$-finitely  presented $A$-module and $B^n$ is a finitely generated  $A$-module (since $B$ is a  finitely generated $A$-module),   $K$ is an $V$-finite $A$-module, and so  $K$ is a  $\phi(V)$-finite $B$-module.  Therefore, $M$ is a $\phi(V)$-finitely presented $B$-module.\cqfd

The following result extends \cite[Theorem 2.1.8 (2)]{Glaz}.

\begin{prop}\label{cor-Sfp-Quot}  Let $I$ be an ideal of $R$ and let $M$ be an $R/I$-module. Assume that $I\cap S= \emptyset   $ so that  $T:=\{s+I\in R/I;\, s\in S\}$ is a  multiplicative subset of $R/I$. Then,
\begin{enumerate}
    \item    $M$ is an $S$-finite  $R$-module if and only if $M$ is a $T$-finite  $R/I$-module.
    \item  If $M$ is an $S$-finitely presented $R$-module, then $M$ is a $T$-finitely presented $R/I$-module.
The converse holds when  $I$ is an $S$-finite ideal of $R$.
\end{enumerate}
\end{prop}
\proof 1. Easy.\\
2. Use the canonical ring surjection $R\longrightarrow R/I$ and   Proposition \ref{pro-chan-phi}.\\
Conversely, if $M$ is a $T$-finitely presented $R/I$-module. Then, there is an exact sequence of $R/I$-modules, and then of $R$-modules  $$0\longrightarrow  K \longrightarrow (R/I)^n \longrightarrow M \longrightarrow 0,$$ where $n>0$ is an integer and $K$ is  a $T$-finite $R/I$-module.   By the first assertion, $K$ is also an $S$-finite $R$-module. And since  $I$ is an $S$-finite ideal of $R$,  $(R/I)^n$ is  an $S$-finitely presented  $R$-module. Therefore, by Theorem \ref{thm-f-p-sq} (4), $M$ is an $S$-finitely presented $R$-module.  \cqfd


 
 
\section{$S$-coherent rings}

Before giving the definition of  $S$-coherent rings, we give, following the calssical case, the definition of  $S$-coherent modules.

\begin{defn} \label{def-S-C-mod}\textnormal{An $R$-module $M$ is said to be $S$-coherent, if it is
finitely generated and every finitely generated submodule of $M$ is $S$-finitely presented.}
\end{defn}

Clearly, every coherent module is $S$-coherent. However, using Proposition \ref{pro-S-C-sq}(1) below, one can show that,   for an $S$-finite ideal  $I$ of $R$ which is not finitely generated, the $R$-module $R/I$ is $S$-coherent but it is not coherent.\medskip

The reason of why we consider  finitely generated submodules rather than $S$-finite submodules is explained in   assertion (4) of Remark \ref{rem-S-C-ring}.\medskip

The following result studies the behavior of $S$-coherence of modules 
in short exact sequences.  It generalizes \cite[Theorem 2.2.1]{Glaz}.

\begin{prop}\label{pro-S-C-sq} Let
$0\longrightarrow
P\stackrel{f}{\longrightarrow}N\stackrel{g}{\longrightarrow}
M\longrightarrow 0$ be an exact sequence of $R$-modules.  The 
following assertions hold:
\begin{enumerate}
\item If $P$ is   $S$-finite and $N$ is   $S$-coherent, then
$M$ is   $S$-coherent.
\item If $M$ and $P$ are $S$-coherent, then so is $N$.
In particular, every finite direct sum of   $S$-coherent modules is $S$-coherent.
\item If  $N$ is $S$-coherent and $P$ is finitely generated, then  $P$ is $S$-coherent.
\end{enumerate}
\end{prop}
\proof
1. It is clear that   $M$ is finitely generated. Let $M'$ be a finitely generated  submodule of
$M$. Then, $f(P)\subseteq g^{-1}(M')$, so there exist  two shorts
exacts sequences of $R$-modules\\ $0\longrightarrow K\longrightarrow
R^n\longrightarrow P\longrightarrow0$ and $0\longrightarrow
K'\longrightarrow R^m\longrightarrow M'\longrightarrow0$, where  
 $n$ and $m$ are two positive integers. Then, by Horseshoe Lemma, we get the following diagram 
$$\xymatrix{&0&0&0&\\0\ar[r]&P\ar[r]\ar[u]&g^{-1}(M')\ar[u]\ar[r]&M'\ar[u]\ar[r]
&0\\0\ar[r]&R^n\ar[u]\ar@{.>}[r]&R^{n+m}\ar@{.>}[u]\ar@{.>}[r]&R^m\ar[u]\ar[r]&0
\\0\ar[r]&K\ar[u]\ar[r]&K''\ar[u]\ar[r]&K'\ar[r]\ar[u]&0\\&0\ar[u]&0\ar[u]&0\ar[u]&}$$
Since   $g^{-1}(M')$ is a finitely generated submodule of the $S$-coherent module
$N$,   $g^{-1}(M')$ is $S$-finitely presented. Then, using  Theorem \ref{thm-f-p-sq} (5), 
$K''$ is $S$-finite, and so $K'$ is $S$-finite. Therefore, $M'$ is $S$-finitely presented.\\ 
2. Clearly $N$ is finitely generated. Let $N'$ be a  finitely generated submodule of  $N$. Consider  the exact
sequence $0\longrightarrow Ker (g_{/N'})  \stackrel{f}{\longrightarrow} N'\stackrel{g}{\longrightarrow}
g(N')\longrightarrow 0$. Then, $g(N')$ is a  finitely generated  submodule of  the $S$-coherent module $M$. Then, $g(N')$  is $S$-finitely presented.  Then, $Ker (g_{/N'})$ is finitely generated by  Theorem \ref{thm-f-p-sq} (5), and since $P$ is  $S$-coherent, $Ker (g_{/N'})$ is $S$-finitely presented. Therefore, by
(2) of  Theorem \ref{thm-f-p-sq}, $N'$ is $S$-finitely presented.\\ 
3. Evident since a submodule of $P$ can be seen as a submodule of $N$.\cqfd

Now we set the definiton of $S$-coherent rings.

\begin{defn}\label{def-S-C-ring} \textnormal{A ring $R$ is called   $S$-coherent, if it is $S$-coherent as an $R$-module;  that is, if every finitely generated ideal of $R$ is $S$-finitely presented.}
\end{defn}

\begin{rem}\label{rem-S-C-ring}
\textnormal{
\begin{enumerate}
\item   Note  that every $S$-Noetherian ring is $S$-coherent. Indeed, this follows from the fact that when $R$ is $S$-Noetherian, every finitely generated free $R$-module is $S$-Noetherian (see the discussion before \cite[Lemma 3]{AD02}). Next, in Example \ref{exm-S-Noe-C}, we give an example  of an $S$-coherent ring which is not  $S$-Noetherian.
    \item Clearly, every coherent ring is  $S$-coherent. The converse is not true in general. As an example of an $S$-coherent ring which is not coherent, we consider the trivial extension $A =\Z\ltimes ( \Z/2\Z) ^{(\N)}$ and the multiplicative set $V =\{(2,0)^n;\,n\in\N\}$. Since  $ (0:(2,0))=0\ltimes M$ is not finitely generated, $T$ is not coherent.  Now, for every ideal $I$ of $A$,    $(2,0)I $   is finitely generated; in fact, $(2,0)I=2J\ltimes 0$, where $J=\{a\in \Z;\,\exists b\in  ( \Z/2\Z) ^{(\N)},(a,b)\in I\}$. Since $J$ is an ideal of $Z$, $J=a\Z$ for some element $a\in \Z$. Then,  $(2,0)I=2J\ltimes 0=(2a,0) A$.  This shows that $A$ is $V$-Noetherian and so $V$-coherent.
      \item It is easy to show that, if $M$ is an $S$-finitely presented $R$-module, then $M_S$ is a  finitely presented $R_S$-module. Thus, if $R$ is a  $S$-coherent ring, $R_S$ is a coherent ring. However, it seems not evident to give a condition so that the converse holds, as done for $S$-Noetherian rings (see \cite[Proposition 2 (f)]{AD02}). In Section 4, we give another $S$-version of coherent rings which can be characterized  in terms of localization.
    \item One would propose for an $S$-version of coherent rings, the following condition ``\textbf{$S$-$C$}: every $S$-finite   ideal of $R$ is $S$-finitely presented". However, if $R$ satisfies the condition \textbf{$S$-$C$}, then in particular, every  $S$-finite   ideal of $R$ is finitely generated. So, every  $S$-finite ideal of $R$ is finitely presented; in particular, $R$ is coherent. This means that the notion of rings  with the  condition \textbf{$S$-$C$} cannot be considered as an $S$-version of the classical coherence. Nevertheless, these rings could be of particular interest as a new class of rings between the class of coherent rings and the class of Noetherian rings.\\ To give an example of a coherent ring which does not satisfy the condition \textbf{$S$-$C$}, one could consider the Boolean ring $B=\prod\limits_{i=1}^{\infty} k_i$, where $k_i$ is the field of two elements for every $i\in \N$, and the  multiplicative subset $V=\{1,e\}$ of $B$, where   $e=(1,0,0,...)\in B$. Indeed, the ideal $B=\bigoplus\limits_{i=1}^{\infty} k_i$ is $V$-finite but not finitely generated.\\
Also, note that the following condition ``\textbf{$S$-$c$}: every $S$-finite   ideal of $R$ is finitely generated" could be of interest. Indeed, clearly one can show the following equivalences:
\begin{enumerate}
    \item A ring $R$ satisfies the condition \textbf{$S$-$C$} if and only if $R$ is coherent and  satisfies the condition \textbf{$S$-$c$}. 
    \item A ring $R$ is coherent   if and only if $R$ is $S$-coherent and  satisfies the condition \textbf{$S$-$c$}.
 \item A ring $R$ is Noetherian   if and only if $R$ is $S$-Noetherian and  satisfies the condition \textbf{$S$-$c$}.
\end{enumerate}  
\end{enumerate}
}
\end{rem}

To give  an example  of an $S$-coherent ring which is not  $S$-Noetherian, we use the following result.

\begin{prop}\label{pro-S-C-sum}  Let $R=\displaystyle\prod_{i=1}^n R_i$ be a direct product of rings $R_i$ ($n\in \N $) and $S=\displaystyle\prod_{i=1}^n S_i$ be a  cartesian product of multiplicative sets $S_i$ of $R_i$. Then, $R$ is $S$-coherent if and only if $R_i$ is $S_i$-coherent for every $i\in\{1,...,n\}$.
\end{prop}
\proof The result is proved using  standard arguments.\cqfd

\begin{exmp}\label{exm-S-Noe-C}
Consider the ring $A$ given in Remark \ref{rem-S-C-ring} (2). Let $B$ be a coherent ring which has a multiplicative set $W$ such that $B_V$ is not Noetherian. Then, $A\times B$ is $V\times W$-coherent (by Proposition \ref{pro-S-C-sum}), but it is not $V\times W$-Noetherian (by \cite[Proposition 2 (f)]{AD02}).\medskip
\end{exmp}

Now, we give our main result. It is the $S$-counterpart of the classical Chase's result \cite[Theorem  2.2]{C60}. We mimic the proof of \cite[Theorem 2.3.2]{Glaz}. So we use the following lemma.

\begin{lem}[\cite{Glaz}, Lemma 2.3.1]\label{lem-princ} Let $R$ be a ring, let $I=(u_1,u_2,...,u_n)$
be a  finitely generated ideal of $R$ ($n\in \N$) and let $a\in R$. Set $J=I+Ra$. Let $F$ be a free module
on generators $x_1,x_2,...,x_{r+1}$ and let $0\longrightarrow
K\longrightarrow F\stackrel{f}{\longrightarrow} J\longrightarrow 0$,
be an exact sequence with $f(x_i)=u_i$ $(1\leq i\leq r)$ and
$f(x_{r+1})=a$. Then there exists an exact sequence $0\longrightarrow
K\cap F'\longrightarrow K\stackrel{g}{\longrightarrow}
(I:a)\longrightarrow 0$, where $F'=\bigoplus\limits_{i=1}^nRx_i$.
\end{lem}

\begin{thm} \label{thm-princ} The following assertions are equivalent:
\begin{enumerate}
\item $R$ is $S$-coherent.
\item Every $S$-finitely presented $R$-module is $S$-coherent.
\item Every finitely generated $R$-submodule of  a free $R$-module is $S$-finitely presented.
\item $(I:a)$ is an $S$-finite ideal of $R$, for every finitely generated ideal $I$ of $R$ and $a\in R$.
\item $(0:a)$ is an $S$-finite ideal of $R$ for every $a\in R$ and the intersection of two finitely generated  ideals of $R$ is an $S$-finite ideal of $R$.
\end{enumerate}
\end{thm}
\proof
The proof is similar to that of \cite[Theorem  2.2]{C60} (see also \cite[Theorem 2.3.2]{Glaz}). However, for the sake of completeness we give its proof here. \\
(1$\Rightarrow$2) Follows from Proposition \ref{pro-S-C-sq} (1).\\ 
(2$\Rightarrow$1) Obvious.\\
(1$\Rightarrow$3) Let $N$ be a finitely generated submodule of a free $R$-module $F$. Hence, there exists  a  finitely generated free submodule $F'$ of $F$ containing $N$. Then, by (1), $F'$ is  $S$-coherent. Therefore, $N$ is $S$-finitely presented.\\
(3$\Rightarrow$1) Trivial.\\
(1$\Rightarrow$4) Let $I$ be a finitely generated  ideal of $R$. Then, $I$ is $S$-finitely presented. Consider $J=I+Ra$, where $a\in R$. Then, $J$   is   finitely generated, and so  it is $S$-finitely presented. Thus,	 there exists an exact sequence  $0\longrightarrow K\longrightarrow R^{n+1}\longrightarrow J\longrightarrow 0$, where
$K$ is $S$-finite. By Lemma \ref{lem-princ}, there exists a surjective homomorphism $g:K\longrightarrow(I:a)$ which shows that $(I:a)$ is $S$-finite.\\
($4\Rightarrow 1$) This is proved by induction on $n$, the   number of   generators of a  finitely generated ideal $I$ of $R$. For $n=1$, use assertion $(4)$ and the exact sequence $0\longrightarrow  (0:I)
 \longrightarrow R \longrightarrow I\longrightarrow 0$. For $n>1$, use assertion $(4)$ and  Lemma \ref{lem-princ}.\\
($1\Rightarrow 5$) Since $R$ is $S$-coherent, Proposition \ref{prop-S-fp-inde} applied on the exact sequence $0\longrightarrow (0:a) \longrightarrow R \longrightarrow aR\longrightarrow 0$ shows that the ideal $(0:a)$ is $S$-finite. Now, Let  $I$ and $J$ be two  finitely generated   ideals of $R$. Then, $I+J$ is  finitely generated and so $S$-finitely presented. Then,   applying Theorem \ref{thm-f-p-sq} (5) on the short the exact sequence  $0\longrightarrow I\cap J\longrightarrow I\oplus J\longrightarrow I+J\longrightarrow
0$, we get that  $I\cap J$ is $S$-finite.\\
($5\Rightarrow 1$) This is proved by induction on the number of   generators of a  finitely generated ideal $I$ of $R$, using the two short exact sequences used in $1\Rightarrow 5$.\cqfd
 
It is worth noting that, in Chase's paper \cite{C60}, coherent rings were characterized using the  notion of flat modules. Then, naturaly one can ask of an $S$-version of flatness that characterizes $S$-coherent rings similarly to the classical case. We leave it as an interesting open question.\medskip

We end this section with some change of rings results.\medskip

The following results extends \cite[Theorem 2.4.1]{Glaz}.

\begin{prop}\label{prop-SC-quotien} Let $I$ be an $S$-finite ideal of $R$. Assume that $I\cap S=\emptyset $ so that  $T:=\{s+I\in R/I;\, s\in S\}$ is a  multiplicative subset of $R/I$. Then, an $R/I$-module $M$ is $T$-coherent if only if it is an  $R$-module $S$-coherent. In particular, the following assertions hold:
\begin{enumerate}
\item If $R$ is an $S$-coherent ring, then $R/I$ is a  $T$-coherent ring.
\item If $R/I$ is a $T$-coherent ring and $I$ is an $S$-coherent $R$-module, then $R$ is an $S$-coherent ring.
\end{enumerate}
\end{prop}
\proof  Use Proposition \ref{cor-Sfp-Quot}.\cqfd

Next result generalizes \cite[Theorem 2.4.2]{Glaz}. It studies the transfer of $S$-coherence under localizations.

\begin{lem} \label{lem-SC-loca}
 Let $f:A\to B $ be a ring homomorphism such that $B$ is a flat $A$-module, and let $V$ be a multiplicative set of $A$. If an $A$-module $M$ is $V$-finite (resp., a $V$-finitely presented), then $M\otimes_A B$ is an  $f(V)$-finite (resp., $f(V)$-finitely presented) $B$-module.
\end{lem}
\proof Follows using the fact that flatness preserves injectivity.\cqfd

\begin{prop}  \label{prop-SC-loca}
If $R$ is $S$-coherent, then  $R_T$ is an $S_T$-coherent ring for every multiplicative set $T$ of $R$.
\end{prop}
\proof Let $J$ be a finitely generated ideal of $R_T$. Then, there is a  finitely generated ideal   $I$ of $R$ such that $J=I_T$. Since  $R$ is $S$-coherent, $I$ is $S$-finitely presented. Then, using Lemma \ref{lem-SC-loca}, the ideal  $J=I \otimes_R R_T$ of $R_T$ is  $S_T$-finitely presented, as desired.\cqfd


\section{Other   $S$-version of  finiteness}

In this short section, we present another   $S$-version of $S$-finiteness and we prove that this notion can be characterized  in terms of localization.\medskip 

The following definition gives  another $S$-version of finitely presented modules.
 
\begin{defn}\label{def-c-S-fp} \textnormal{An $R$ module $M$ is called  c-$S$-finitely presented, if there exists a finitely presented submodule $N$ of $M$  such that $sM\subseteq N  \subseteq M$ for some   $s\in S$.}
\end{defn}

\begin{rem}\label{rem-def-c-S-fp} 
\begin{enumerate}
    \item  Clearly, every finitely presented module is c-$S$-finitely presented. However, the converse does not hold in general. For that it suffices to consider a coherent ring which has an $S$-finite module which is not finitely generated.  An example of a such ring is given in Remark \ref{rem-S-C-ring} (4).  
 \item  The inclusions in Definition \ref{def-c-S-fp}  complicate  the study of the behavior of of c-$S$-finitely presented modules  in short exact sequences as done in Theorem \ref{thm-f-p-sq}. This is why we think that c-$S$-finitely presented modules will be mostly used by commutative rings theorists rather than researchers interested in notions of homological algebra. This is the reason behind the use of the letter ``c" in ``c-$S$-finitely presented".
  \item  It seems that there is not any relation between the two notions of c-$S$-finitely presented   and   $S$-finitely presented modules. Nevertheless, we can deduce that in a c-$S$-coherent ring (defined below), every $S$-finitely presented ideal  is c-$S$-finitely presented.  
\end{enumerate}
\end{rem}

It is well-known that if, for an $R$-module $M$, $M_S$ is a finitely presented $R_S$-module, then there is a finitely presented $R$-module $N$ such that $M_S=N_S$.  Nevertheless, what doest not  make things work with respect to localization for $S$-finitely presented modules is the fact that the module $N$ which satisfies $M_S=N_S$ is not necessarily a submodule of $M$. For c-$S$-finitely presented modules we give the following result.

\begin{prop}\label{prop-c-Sfp-loc} 
\begin{enumerate}
    \item If an $R$-module $M$  is  c-$S$-finitely presented, then $M_S$ is a finitely presented $R_S$-module. 
    \item  A  finitely generated $R$-module $M$ is  c-$S$-finitely presented if and only if there is a finitely presented submodule $N$ of $M$  such that $M_S=N_S$.
\end{enumerate}
\end{prop}  
\proof  1. Obvious.\\
2. ($\Rightarrow $) Clear.\\
($\Leftarrow $) Since $M$ is   finitely   and $M_S=N_S$, there is an $s\in S$ such that $sM\subseteq N$, as desired.\cqfd

Now we define the other $S$-version of the classical coherence of rings.

\begin{defn}\label{def-c-S-C-ring} \textnormal{A ring $R$ is called   c-$S$-coherent, if  every $S$-finite  ideal of $R$ is $S$-finitely presented.}
\end{defn}

Clearly, every coherent ring is  c-$S$-coherent. The converse is not true in general. The ring given in Example \ref{rem-S-C-ring} (2) can be used as an example of a c-$S$-coherent ring which is not coherent.\medskip

Also, it is evident that every  $S$-Noetherian ring is  c-$S$-coherent. As done in Example \ref{exm-S-Noe-C}, we use the following result to give  an example  of  a c-$S$-coherent  ring which is not   $S$-Noetherian.

\begin{prop}\label{pro-c-S-C-sum}  Let $R=\displaystyle\prod_{i=1}^n R_i$ be a direct product of rings $R_i$ ($n\in \N $) and $S=\displaystyle\prod_{i=1}^n S_i$ be a  cartesian product of multiplicative sets $S_i$ of $R_i$. Then, $R$ is c-$S$-coherent if and only if $R_i$ is c-$S_i$-coherent for every $i\in\{1,...,n\}$.
\end{prop}
\proof The result is proved using  standard arguments.\cqfd

\begin{exmp}\label{exm-S-Noe-c-C}
Consider a c-$S$-coherent ring $A$ which is not coherent. Let $B$ be a coherent ring which has a multiplicative set $W$ such that $B_V$ is not Noetherian. Then, $A\times B$ is c-$V\times W$-coherent (by Proposition \ref{pro-c-S-C-sum}), but it is not $V\times W$-Noetherian (by \cite[Proposition 2 (f)]{AD02}).\medskip
\end{exmp}

The follwoing result characterizes c-$S$-coherent rings can be characterized  in terms of localization.

\begin{thm} \label{thm-princ2} The following assertions are equivalent:
\begin{enumerate}
\item $R$ is c-$S$-coherent.
\item Every  finitely generated ideal of  $R$  is  c-$S$-finitely presented.
\item For every finitely generated ideal $I$ of $R$, there is a finitely presented ideal $J\subseteq I$ such that $I_S=J_S$. In particular,  $R_S$ is a coherent ring.
\end{enumerate}
\end{thm}
\proof
(1$\Rightarrow2 \Rightarrow 3$ ) Straightforward.\\ 
(3$\Rightarrow$1) Let $I$ be an $S$-finite ideal of $R$. Then,  there exist an $s\in S$ and a finitely generated ideal  $J$ of  $R$  such that $sI\subseteq J\subseteq I$. By assertion $(3)$, there is a finitely presented ideal $K\subseteq J$ such that $K_S=J_S$. Then, there is a $t\in S$ such that $tJ\subseteq K$. Therefore, $tsI \subseteq    K \subseteq I$, as desired.\cqfd

We end the paper with a result which relates c-$S$-coherent rings with the notion of $S$-saturation.\medskip

 In \cite{AD02}, the notion of $S$-saturation is used to characterize $S$-Noetherian rings. Assume that $R$ is an integral domain.  Let $Sat_S(I)$ denotes the $S$-saturation of an ideal $I$ of $R$; that is, $Sat_S(I):=IR_S\cap R$. In \cite[Proposition 2 (b)]{AD02}, it is proved that if $Sat_S(I)$ is $S$-finite, then $I$ is $S$-finite and $Sat_S(I)=(I:s)$ for some $s\in S$.  This fact was used to prove that a ring $R$ is $S$-Noetherian if and only if $R_S$ is Noetherian and, for every finitely generated ideal of $R$,  $Sat_S(I)=(I:s)$  for some $s\in S$ (see \cite[Proposition 2 (f)]{AD02}). The following result shows that the implication of \cite[Proposition 2 (b)]{AD02} is in fact an equivalence in more general context. \medskip

Consider $ N\subseteq M$ an inclusion  of  $R$-modules. Let $f:M \to M_S$ be the canonical $R$-module homomorphism. Denote by $f(N)R_S$ the $R_S$-submodule of $M_S$ generated by  $f(N)$. We set  $Sat_{S,M}(N):= f^{-1} (f(N)R_S)$ and $(N:_M s):=\{m\in M; \, sm\in N\}$.

\begin{prop}\label{prop-Sat-S-f} Let $N$ be an $R$-submodule of  an $R$-module  $M$.   $Sat_{S,M}(N)$ is $S$-finite if and only if   $N$ is $S$-finite and $Sat_{S,M}(N)=(N:_M s)$ for some $s\in S$. 
\end{prop} 
\proof ($\Rightarrow$) Set $K=Sat_{S,M}(N)$. Since $K$ is $S$-finite, there exist  an $s\in S$ and a finitely generated $R$-module   $J$ such that  $sK\subseteq J\subseteq K$. Thus, $sN\subseteq sK\subseteq J$.  We can write $J=Rx_1+Rx_2+\cdots+Rx_n$ for some $x_1,x_2,...,x_n\in J$. For each $x_i$, there exists a $t_i\in S$ such
that $t_ix_i\in N$. We set $t=\prod\limits_{i=1}^nt_i$. Then,  $ts N\subseteq ts K\subseteq tJ\subseteq N$. Then, $N$ is $S$-finite. On the other hand,  since $s K\subseteq tJ \subseteq N\subseteq K$,   $K\subseteq(N:_M s)$. Conversely, let $x\in(N:_M s)$. Then, $sx\in N$, so  $x\in K$, as desired.\\
($\Leftarrow$)  Since $N$ is $S$-finite, there exist a $t\in S$ and a finitely generated $R$-module  $J$ such that $tN\subseteq J\subseteq N$. On the other hand, since   $K=(N:s)$ for some $s\in S$, $s K\subseteq N$.  Consequently, $ts K\subseteq tN\subseteq J\subseteq N\subseteq K$. Therefore, $K$ is $S$-finite.\cqfd

The following result is proved similarly to the  proof of Proposition \ref{prop-Sat-S-f}. However, to guarantee the preservation of finitely presented modules when multiplying by elements of $S$, we assume that $S$ does not contain any zero-divisor of $R$.

\begin{prop}\label{prop-Sat-S-fp} Assume that every element of $S$ is regular. Let $N$ be an $R$-submodule of an $R$-module  $M$. $Sat_{S,M}(N)$ is c-$S$-finitely presented if and only if $N$ is c-$S$-finitely presented and $Sat_{S,M}(N)=(N:_M s)$ for some $s\in S$. 
\end{prop}

\begin{cor}\label{cor-princ2} Assume that every element of $S$ is regular. 
The following assertions are equivalent:
\begin{enumerate}
\item For every  finitely generated ideal $I$ of $R$, $Sat_S(I)$ is  c-$S$-finitely presented.
\item $R$ is c-$S$-coherent and, for every finitely generated ideal $I$ of $R$,   $Sat_S(I)=(I:s)$  for some $s\in S$.
\end{enumerate}
\end{cor}

\medskip


 \noindent {\bf Acknowledgement.} A part of this work was presented by the second author at the Scientific day of Algebra ``JA GRAAF 2016" held in Faculty of Sciences of Rabat (May 17, 2016).\\
The authors would like  to thank Professor  Zine El Abidine Abdelali  for his helpful comments during the preparation of this paper.
\medskip

\bibliographystyle{amsplain}

\end{document}